\documentclass[12pt]{article}
\usepackage{amsmath,amssymb,amsthm, amscd}
\usepackage[dvipdfmx]{graphicx}
\usepackage{color}
\usepackage{here}

\usepackage{mathrsfs}
  
  \makeatletter
  \@addtoreset{equation}{section}
  \makeatother

\renewcommand{\tilde}[1]{\widetilde{#1}}
\newtheorem{dfn}{Definition}[section]
\newtheorem{thm}[dfn]{Theorem}
\newtheorem{lem}[dfn]{Lemma}
\newtheorem{prop}[dfn]{Proposition}

\begin{document}
\begin{center} 
{\bf {\LARGE Singularities of Fano varieties of lines on singular cubic fourfolds}}
\end{center} 
\vspace{0.4cm}

\begin{center}
{\large Ryo Yamagishi}
\end{center} 
\vspace{0.4cm}

\begin{abstract}
Let $X$ be a cubic fourfold that has only simple singularities and does not contain a plane. We prove that the Fano variety of lines on $X$ has the same analytic type of singularity as the Hilbert scheme of two points on a surface with only ADE-singularities. This is shown as a corollary to the characterization of a singularity that is obtained as a $K3^{[2]}$-type contraction and has a unique symplectic resolution.
\end{abstract}

\section{Introduction}\label{1}

When one wants to construct irreducible holomorphic symplectic manifolds, one method is to consider moduli spaces of subschemes, sheaves or more generally objects in the derived categories of smooth cubic fourfolds. For example, Beauville and Donagi showed that the Fano variety of lines  on a smooth cubic fourfold is an irreducible holomorphic symplectic manifold \cite{BD}. Other examples are given in \cite{LLSvS}, \cite{LMS}, \cite{Ou}, and so on. Then it is natural to consider singular cubic fourfolds in order to obtain singular symplectic varieties. Indeed, C. Lehn showed in \cite{Le} that the Fano variety $F(X)$ of lines on a cubic fourfold $X$ that has only simple singularities and does not contain a plane has symplectic singularities in the sense of \cite{B} (see \S2). In this article we will study the local structure of $F(X)$ and show that it has the same singularity type as the Hilbert scheme $\mathrm{Hilb}^2(\Gamma)$ of 2 points on a surface $\Gamma$ with only ADE-singularities (Theorem \ref{thm:main1} and \ref{thm:main2}). As a corollary, we see that the germs of singular points in $F(X)$ only depends on the ADE-types of the singular points of $X$. This seems interesting since even local structure of $F(X)$ is a priori determined by the global structure of $X$.

To state the main theorem more precisely, we should note that (singular) symplectic varieties have relatively simple structures. Kaledin showed that any symplectic variety has natural stratification by locally closed symplectic submanifolds \cite[Thm. 2.3]{K}. In our case $F(X)$ is 4-dimensional, and the 2-dimensional singular locus of $F(X)$ is divided into finitely many ``bad'' points $q_1,\dots,q_k$ and the complementary open subset $U$. The theorem of Kaledin also tells us that the local structure of $F(X)$ at each point in $U$ is the same as the product $\mathbb{C}^2\times\text{(ADE-singularity)}$. Therefore, the problem is to study the structure near $q_i$'s. In \cite{Y}, it is shown that $\mathrm{Hilb}^2(\Gamma)$ is also a symplectic variety and thus it has stratification as well. We will show that for each $q_i\in F(X)$ there is some ``bad'' point $q$ in the singular locus of $\mathrm{Hilb}^2(\Gamma)$ such that the local structures near $q_i$ and $q$ are the same, i.e., the formal completions at these points are isomorphic.

The main theorems will be shown as a corollary to a slightly more general result (Theorem \ref{thm:main3}). This theorem says that a symplectic variety $Z$ whose unique symplectic resolution is a symplectic manifold of $K3^{[2]}$-type has the same singularity type as $\mathrm{Hilb}^2(\Gamma)$.

To prove this theorem, we use the Bayer-Macr\`{i} theory about moduli spaces of Bridgeland stable objects on a twisted $K3$ surface (cf. \cite{BM}). Then we can determine the structure of a special fiber of a symplectic resolution of $Z$ and show that it is isomorphic to a special fiber for $\mathrm{Hilb}^2(\Gamma)$. By using the result in \cite{Y}, we deduce that $Z$ and $\mathrm{Hilb}^2(\Gamma)$ have the same singularities. 

\section{Generalities on $F(X)$}\label{2}

Throughout this article all schemes are defined over the complex number field $\mathbb{C}$. By the {\em singularity type} (or the {\em analytic type})  of a variety at a point, we mean the isomorphism class of the formal completion of the variety at the point.

Let $X\subset \mathbb{P}^5$ be a cubic hypersurface. In this section we assume that $X$ has only simple singularities in the sense of Arnold \cite{Ar}, i.e., every singular point $p$ of $X$ is an isolated singularity defined by one of the following equations
$$\begin{aligned}
\mathbf{A}_n\;&:\;x_1^2+x_2^2+x_3^2+x_4^2+x_5^{n+1}\;(n\ge1)\\
\mathbf{D}_n\;&:\;x_1^2+x_2^2+x_3^2+x_4^2x_5+x_5^{n-1}\;(n\ge4)\\
\mathbf{E}_6\;&:\;x_1^2+x_2^2+x_3^2+x_4^3+x_5^4\\
\mathbf{E}_7\;&:\;x_1^2+x_2^2+x_3^2+x_4^3+x_4x_5^3\\
\mathbf{E}_8\;&:\;x_1^2+x_2^2+x_3^2+x_4^3+x_5^5\\
\end{aligned}$$
for some analytically local coordinates $x_1,x_2,\dots,x_5$ near $p$.

We fix a singular point $p\in X$ and assume that $p=[0\mathbin:0\mathbin:0\mathbin:0\mathbin:0\mathbin:1]\in\mathbb{P}^5$ by a linear change of coordinates.  Let $x_0,\dots,x_5$ be the homogeneous coordinates of $\mathbb{P}^5$. By the condition that $p$ is a singular point of $X$, one sees that $X$ is defined by the following cubic polynomial
\begin{eqnarray}
f(x_0,\dots,x_5)=x_5f_2(x_0,\dots,x_4)+f_3(x_0,\dots,x_4)
\label{eq1}
\end{eqnarray}
where $f_2(x_0,\dots,x_4)$ and $f_3(x_0,\dots,x_4)$ are homogeneous polynomials of degree 2 and 3 respectively. 

Let $F(X)=\{l\in \mathrm{Grass}(1,\mathbb{P}^5) \mid l\subset X\}$ be the Fano variety of lines on $X$. $F(X)$ contains a subvariety $S_p=\{l\in \mathrm{Grass}(1,\mathbb{P}^5) \mid p\in l\subset X\}$. The natural map from $S_p$ to the hyperplane $H=\{x_5=0\}\subset\mathbb{P}^5$ given by $l\mapsto l\cap H$ is an embedding and one can check that the image in $H\cong\mathbb{P}^4$ with the coordinates $x_0,\dots,x_4$ is defined by $f_2=f_3=0$. 

From now on we assume that $X$ contains no planes and in particular $S_p$ contains no lines. It is shown in \cite[Prop. 5.8(b)]{O'G} and \cite[Lem. 3.3]{Le} that $S_p$ is a normal surface with only ADE-singularities and that its minimal resolution is a $K3$ surface. We can determine the singularity type of $S_p$ from that of $X$ by the result of Wall. To state the result, one should note that the exceptional divisor $E$ of the blowing-up $\tilde{X}$ of $X$ at $p$ is isomorphic to the quadric $\mathcal{Q}=\{f_2=0\}$ in $\mathbb{P}^4$. We put $\mathcal{C}=\{f_3=0\}\subset\mathbb{P}^4$.

\begin{prop}(Wall, \cite[Thm. 2.1]{Wa})\label{prop:Wall}
Let $q\in S_p=\mathcal{Q}\cap\mathcal{C}$ be a singular point. Then the followings hold.\\
(1)\,Either ``$q\in \mathrm{Sing}(\mathcal{Q})$ and $q\notin \mathrm{Sing}(\mathcal{C})$'' or ``$q\notin \mathrm{Sing}(\mathcal{Q})$'' holds.\\
(2)\,If $q\in \mathrm{Sing}(\mathcal{Q})$, then there does not exist a singular point of $X$ on the line $\overline{pq}$ (here $q$ is regarded as a point of $X$) other than $p$ and $q$ and the ADE-type of $S_p$ at $q$ is the same as that of $\tilde{X}$ at the point of $E$ that corresponds to $q$.\\
(3)\,If $q\notin \mathrm{Sing}(\mathcal{Q})$, then there exists exactly one singular point $p'$ of $X$ on $\overline{pq}$ other than $p$ and $q$ and the ADE-type of $S_p$ at $q$ is the same as that of $X$ at $p'$.
\end{prop}

By this proposition we can explicitly describe singular points of $S_p$. Let $\{p,p_1,\dots,p_k\}$ ($k$ can be 0) be the set of all singular points of $X$ and assume that  the singularity types of these points are $T_n,T_{n_1},\dots,T_{n_k}$ respectively where $T=\mathbf{A},\mathbf{D}$ or $\mathbf{E}$. From the equations of simple singularities one knows that the corank of the quadric $\mathcal{Q}$ and the types of singular points of $\tilde{X}$ on $E$ are determined as in Table \ref{table1}. In the table one should understand $\mathbf{A}_0$ and $\mathbf{D}_3$ as $\emptyset$ and $\mathbf{A}_3$ respectively.\\

\begin{table}[H]
\begin{center}
\caption{Degeneracy of quadrics and singularities of $\tilde{X}$}\label{table1}
\begin{tabular}{|c||c|c|c|}
\hline
T&corank of $\mathcal{Q}$&$\mathrm{Sing}(E)$&$\mathrm{Sing}(\tilde{X})\cap E$\\
\hline\hline
$\mathbf{A}_1$&0\,(nondegenerate)&$\emptyset$&$\emptyset$\\ \hline
$\mathbf{A}_n\,(n\ge 2)$&1&point&$\mathbf{A}_{n-2}$\\
\hline
$\mathbf{D}_4$&2&line&$3\mathbf{A}_1$\\
\hline
$\mathbf{D}_n\,(n\ge 5)$&2&line&$\mathbf{A}_1+\mathbf{D}_{n-2}$\\
\hline
$\mathbf{E}_6$&2&line&$\mathbf{A}_5$\\
\hline
$\mathbf{E}_7$&2&line&$\mathbf{D}_6$\\
\hline
$\mathbf{E}_8$&2&line&$\mathbf{E}_7$\\
\hline\end{tabular}

\end{center}
\end{table}

Note also that any line passing through two singular points of $X$ is always contained in $X$ since such a line intersects $X$ with multiplicity at least 4 but $X$ is cubic and the line must be in $X$ by B\'ezout's theorem. By combining these with Proposition \ref{prop:Wall}, one can conclude that $S_p$ has $k$ singular points of types $T_{n_1},\dots,T_{n_k}$ coming from intersections with other $S_{p_i}$'s and extra singular point(s) coming from ``the blowing-up of $T_n$.'' For example, if $T_n=\mathbf{D}_4$, then $S_p$ has 3 more singular points of type $\mathbf{A}_1$ other than the $k$ singular points of types $T_{n_1},\dots,T_{n_k}$.

Next we introduce the results by C. Lehn about $F(X)$. An important property of $F(X)$ is that it is a symplectic variety in the sense of \cite{B}.

\begin{thm}(C. Lehn, \cite[Thm. 3.6]{Le})\label{thm:Fano}
$F(X)$ is a singular symplectic variety birational to $\mathrm{Hilb}^2(S_p)$, the Hilbert scheme of 2 points on $S_p$.
\end{thm}

In the proof of the above theorem in \cite{Le}, a birational map $\mathrm{Hilb}^2(S_p)\dashrightarrow F(X)$ is explicitly constructed as follows. For an element $\xi\in\mathrm{Hilb}^2(S_p)$, the intersection of $X$ and the linear span of $p$ and $\xi$ (regarded as points in the hyperplane $H\subset\mathbb{P}^5$) is the union of 3 lines since $X$ is cubic and $S_p$ does not contain a line by assumption. Two of the lines are the cone over $\xi$ and the residual line gives an element of $F(X)$. Note that this map is well-defined even if $\xi$ is supported at one point. Therefore, we have a birational morphism $\pi:\mathrm{Hilb}^2(S_p)\to F(X)$ and one can check that the indeterminacy locus of the rational map $\pi^{-1}$ is $S_p$. 

We can describe the singular locus of $F(X)$ by the following lemma.

\begin{lem}(cf. \cite[Prop. 3.5]{Le})\label{S_p}
The singular locus of $F(X)$ is $\bigcup_{p\in \mathrm{Sing}(X)} S_p$. Any two of the irreducible components $S_p$ intersect at one point but any three of them do not intersect at one point.
\end{lem}
{\em Proof.} We first show that any line $l\in S_p$ for $p\in \mathrm{Sing}(X)$ is a singular point of $F(X)$. Since $F(X)$ is 4-dimensional, it suffices to show that $\dim T_l F(X)>4$. By a linear change of coordinates, we may assume that $l$ is defined as $\{x_0=x_1=x_2=x_3=0\}$ and $p=[0\mathbin:0\mathbin:0\mathbin:0\mathbin:0\mathbin:1]$ in $\mathbb{P}^5$. Let $f\in H^0(\mathcal{O}_{\mathbb{P}^5}(3))$ be the defining equation of $X$. Then there exist unique $g_i\in H^0(\mathcal{O}_l(2)),\,i=0,1,2,3$ satisfying
$$f=\sum_{i=0}^3 x_ig_i+h$$
where $h\in (x_0,x_1,x_2,x_3)^2$. We have the following exact sequence (see \cite[Ch. 6]{EH}):
$$0\to \mathcal{N}_{l/X}\to \mathcal{O}_l^4(1)\stackrel{\theta}{\to}\mathcal{O}_l(3)$$
where $\theta=(g_0,g_1,g_2,g_3)$.
Since $p$ is a singular point, each $g_i$ is in $x_4\cdot H^0(\mathcal{O}_l(1))$ and thus $H^0(\theta)$ is not surjective. Therefore, $\dim T_l F(X)=\dim H^0(\mathcal{N}_{l/X})>\dim H^0(\mathcal{O}_l^4(1))-\dim H^0(\mathcal{O}_l(3))=8-4=4$.

By the construction of $\pi:\mathrm{Hilb}^2(S_p)\to F(X)$, one sees that $U=\pi^{-1}(F(X)\setminus \bigcup_{p\in \mathrm{Sing}(X)} S_p)$ consists of elements of $\mathrm{Hilb}^2(S_p)$ whose support does not intersect with singular points in $S_p$. Since the Hilbert scheme of points on a smooth surface is again smooth \cite[Thm. 2.4]{Fo}, $U$ and hence $\pi(U)$ are also smooth.

The last claims follow from the fact that the line passing through two singular points of $X$ is contained $X$ and that no three singular points of $X$ are on the same line (see Proposition \ref{prop:Wall} and the argument in the paragraph after it). \qed

\vspace{3mm}

One purpose of this article is to investigate the analytic germs of singular points of $F(X)$. In general Kaledin showed that taking singular loci of symplectic varieties gives stratification by symplectic submanifolds which we call symplectic leaves \cite[Thm. 2.3]{K}. The analytic type of a point $x$ in the 2-dimensional symplectic leaf (or, in other words, the smooth part of $\mathrm{Sing}(F(X))$) is just the same as the product of $\mathbb{C}^2$ and an ADE-singularity. However, when $x\in\mathrm{Sing}(S_p)$, the situation near $x$ in $F(X)$ is not so obvious. The main theorems concern the singularity type near $x$. If $x$ is the intersection point of $S_p$ and another $S_{p_i}$, the singularity type near $x$ is in fact the product of ADE-singularities. This is formulated as follows.

\begin{thm}\label{thm:main1}
Let $X\subset\mathbb{P}^5$ be a cubic hypersurface not containing a plane. Assume $X$ has two simple singular points $p_1$ and $p_2$ of type $T_{n_1}$ and $T_{n_2}$ respectively where $T=\mathbf{A},\,\mathbf{D}$ or $\mathbf{E}$. Set $\{x\}=S_{p_1}\cap S_{p_2}\subset F(X)$. Then the singularity type of $F(X)$ at $x$ is the same as that of $\Gamma_1\times \Gamma_2$ at the unique 0-dimensional symplectic leaf where $\Gamma_i\,(i=1,2)$ is a surface that has a unique singular point of type $T_{n_i}$.
\end{thm}

This will be proved in the next section. If $x$ is not the intersection point, the singularity is a little more complicated. In this case the singularity of $x\in F(X)$ is the same as that of a certain point in the Hilbert scheme of 2 points on a singular surface. As we saw above, $x$ comes from ``the blowing-up of some $T_n$.'' Let $\Gamma$ be a surface that has a unique singular point of type $T_n$. Then the singular locus of the Hilbert scheme $\mathrm{Hilb}^2(\Gamma)$ is isomorphic to the blowing-up of $\Gamma$ at the singular point (cf. \cite[\S2]{Y}). Thus $\mathrm{Hilb}^2(\Gamma)$ has a unique 0-dimensional symplectic leaf if $T=\mathbf{A}$ or $\mathbf{E}$. We denote the singularity type of this leaf by $\bar{\mathbf{A}}_n\,(n\ge3)$ or $\bar{\mathbf{E}}_n\,(n=6,7,8)$ respectively. When $T_n=\mathbf{D}_4$, there are three 0-dimensional symplectic leaves which have the same singularity types in $\mathrm{Hilb}^2(\Gamma)$. We denote it by $\bar{\mathbf{D}}_{4,\mathrm{I}}$. When $T_n=\mathbf{D}_n$ with $n\ge5$, there are two 0-dimensional symplectic leaves. One is of type $\mathbf{A}_1$ in $\mathrm{Sing}(\mathrm{Hilb}^2(\Gamma))$ and the other of type $\mathbf{D}_{n-2}$. We denote the singularity types of them in $\mathrm{Hilb}^2(\Gamma)$ by $\bar{\mathbf{D}}_{n,\mathrm{I}}$ and $\bar{\mathbf{D}}_{n,\mathrm{I\hspace{-1pt}I}}$ respectively. Thus the singularity types which can appear in $\mathrm{Hilb}^2(\Gamma)$ are listed as follows.
\begin{equation}\label{type}
\bar{\mathbf{A}}_n\,(n\ge3),\,\bar{\mathbf{D}}_{n,\mathrm{I}}\,(n\ge4),\,\bar{\mathbf{D}}_{n,\mathrm{I\hspace{-1pt}I}}\,(n\ge5),\text{ and }\bar{\mathbf{E}}_n\,(n=6,7,8).
\end{equation}
The claim about the singularity type of $x$ is stated as follows.

\begin{thm}\label{thm:main2}
Let $X\subset\mathbb{P}^5$ be a cubic hypersurface not containing a plane. Assume $X$ has a simple singular point $p$ of type $T_n$ where $T=\mathbf{A},\,\mathbf{D}$ or $\mathbf{E}$. If $x\in S_p$ is a 0-dimensional symplectic leaf that is not contained in any other irreducible component of $\mathrm{Sing}(F(X))$,  then we have
$$(F(X),x)=
\begin{cases}
\bar{\mathbf{A}}_n&\text{ if}\hspace{5mm}T=\mathbf{A}\\
\bar{\mathbf{D}}_{n,\mathrm{I}}\text{ or }\bar{\mathbf{D}}_{n,\mathrm{I\hspace{-1pt}I}}&\text{ if}\hspace{5mm}T=\mathbf{D}\\
\bar{\mathbf{E}}_n&\text{ if}\hspace{5mm}T=\mathbf{E}.
\end{cases}$$
where being whether $\bar{\mathbf{D}}_{n,\mathrm{I}}$ or $\bar{\mathbf{D}}_{n,\mathrm{I\hspace{-1pt}I}}$ depends on whether the germ $(S_p,x)$ is of type $\mathbf{A}_1$ or $\mathbf{D}_{n-2}$ respectively.
\end{thm}

This will also be proved in the next section.

\section{Proof of the main theorems}\label{3}

In this section we will prove Theorem \ref{thm:main1} and \ref{thm:main2}. We consider a slightly more general situation. Let $Y$ be a projective irreducible holomorphic symplectic manifold of $K3^{[2]}$-type, i.e. deformation equivalent to the Hilbert scheme of two points on a $K3$ surface (see e.g. \cite[Part III]{GHJ}). We assume that we are given a projective birational morphism $\pi:Y\to Z$ where $Z$ is a normal variety and $\pi_*\mathcal{O}_Y=\mathcal{O}_Z$. Then $Z$ is a symplectic variety by definition. In addition, we assume $\pi$ is a unique symplectic resolution (i.e. a resolution by a symplectic manifold) of $Z$.

We will deduce Theorem \ref{thm:main1} and \ref{thm:main2} as a corollary to the following.

\begin{thm}\label{thm:main3}
Assume we are given $\pi:Y\to Z$ as above. If a fiber $\pi^{-1}(z)$ of $z\in Z$ is two-dimensional, then the germ $(Z,z)$ is isomorphic to $(\mathrm{Hilb}^2(\Gamma),q)$ for some $q\in\mathrm{Hilb}^2(\Gamma)$ where $\Gamma$ is a surface with at worst (possibly more than one) ADE-singularities.
\end{thm}

The proof of this theorem will be done in this section by taking several steps. The important observation is that $Y$ is a moduli space of Bridgeland stable objects in the derived category of a (twisted) $K3$ surface and thus its contractions are described in the framework of the Bayer-Macr\`{i} theory \cite{BM}.

Let us recall briefly the theory of Bridgeland moduli spaces over $K3$ surfaces. Refer to \cite[\S2]{BM} for details. Let $S$ be an algebraic $K3$ surface and $\alpha\in \mathrm{Br}(S)$ be a Brauer class. The pair $(S,\alpha)$ is called a twisted $K3$ surface. The cohomology group
$$H^*(S,\mathbb{Z})=H^0(S,\mathbb{Z})\oplus H^2(S,\mathbb{Z})\oplus H^4(S,\mathbb{Z})$$
admits a weight-2 Hodge structure (which is the usual one on $S$ when $\alpha$ is trivial) and we denote the integral $(1,1)$-part by $H_\mathrm{alg}^*(S,\alpha,\mathbb{Z})$. Every object $A$ of $D^b(S,\alpha)$ defines its Mukai vector $v(A)\in H_\mathrm{alg}^*(S,\alpha,\mathbb{Z})$. $H_\mathrm{alg}^*(S,\alpha,\mathbb{Z})$ admits a symmetric bilinear form called the Mukai pairing. With this pairing, $H_\mathrm{alg}^*(S,\alpha,\mathbb{Z})$ becomes an even lattice of signature $(2,\rho)$ where $\rho$ is the Picard number of $S$.

For a generic stability condition $\sigma\in \mathrm{Stab}^\dagger(S)$ and a primitive element $v$ of $H_\mathrm{alg}^*(S,\alpha,\mathbb{Z})$, we can construct a moduli space $M_\sigma(v)$ of $\sigma$-stable objects in $D^b(S,\alpha)$ whose Mukai vectors are $v$. If $v^2\ge -2$, then $M_\sigma(v)$ is of dimension $2n=v^2+2$ and is a projective holomorphic symplectic manifold of $K3^{[n]}$-type. The N\'{e}ron-Severi group $NS(M_\sigma(v))$ with the Beauville-Bogomorov form is identified with $v^\perp\subset H_\mathrm{alg}^*(S,\alpha,\mathbb{Z})$ as a lattice.

In general, a $K3^{[n]}$-type manifold $M$ with $n\ge2$ has a natural extension of the lattice and Hodge structure $H^2(M,\mathbb{Z})\subset \tilde{\Lambda}(M)$ such that the orthogonal complement $H^2(M,\mathbb{Z})^\perp$ is generated by a primitive element $v'$ which is in the $(1,1)$-part with $v^2=\dim M-2$ \cite[\S9]{Ma1}. It is also known that two $K3^{[n]}$-type manifolds $M$ and $M'$ are birational if and only if there is a Hodge isometry of $\tilde{\Lambda}(M)$ and $\tilde{\Lambda}(M')$ which maps $H^2(M,\mathbb{Z})$ isomorphically to $H^2(M',\mathbb{Z})$. If $M$ is a moduli space of Bridgeland stable objects on $(S,\alpha)$, then we have an isomorphism $H^*(S,\alpha,\mathbb{Z})\cong \tilde{\Lambda}(M)$ which sends the Mukai vector to $v'$. 

Now we return to the situation in Theorem \ref{thm:main3}.

\begin{lem}\label{lem:moduli}
$Y$ is a moduli space of Bridgeland stable objects on a twisted $K3$ surface $(S,\alpha)$ for some Mukai vector $v\in H_\mathrm{alg}^*(S,\alpha,\mathbb{Z})$ with $v^2=2$ (and for a suitable stability condition).
\end{lem}

{\em Proof.} By \cite[Thm. 1.2(c)]{BM}, it suffices to show that $Y$ is birational to a moduli space of Bridgeland stable objects. Let $\tilde{\Lambda}_\mathrm{alg}(M)$ be the integral $(1,1)$-part of $\tilde{\Lambda}(M)$. By \cite[Lem. 2.6]{H1} and the argument of the proof of \cite[Prop. 4.1]{H1} (see also \cite[Prop. 4]{Ad}), we only have to show that $\tilde{\Lambda}_\mathrm{alg}(Y)$ contains a sublattice that is isomorphic to a nonzero multiple of the hyperbolic plane $U=\begin{pmatrix}0&1\\1&0\end{pmatrix}$.

We first claim that $\pi$ is a divisorial contraction. Otherwise $\pi$ would contract a $\mathbb{P}^2$ \cite[Thm. 1.1]{WW}. However, this is contrary to the uniqueness assumption of symplectic resolutions since any $\mathbb{P}^2$ can be flopped to give a different symplectic resolution.

Let $\delta'\in NS(Y)\subset\tilde{\Lambda}_\mathrm{alg}(Y)$ be the class of a prime $\pi$-exceptional divisor $E$. Note that $\delta'$ is orthogonal to the pullbacks of the divisors of $Z$. Since the signature of the Beauville-Bogomorov form on $NS(Y)$ is $(1,\rho(Y)-1)$ and the pullback $H$ of an ample divisor of $Z$ satisfies $H^2>0$, we have $\delta'^2<0$. Write $\delta'=m\delta$ with $m\in\mathbb{N}$ and primitive $\delta\in\tilde{\Lambda}_\mathrm{alg}(Y)$. Then $\delta^2=-2$ by \cite[Thm. 1.2]{Ma2}. We can take a primitive generator $v'$ of $H^2(Y,\mathbb{Z})^\perp$ in $\tilde{\Lambda}(Y)$ as mentioned above. Then $v'+\delta$ and $v'-\delta$ generate a lattice which is isomorphic to $U(4)$.
\qed

\vspace{3mm}

Let $E_1,\dots,E_n\subset Y$ be the irreducible exceptional divisors of $\pi$. 

\begin{lem}
The relative Picard number of $\pi$ is equal to $n$ ($=$the number of the exceptional divisors).
\end{lem}

{\em Proof.} Since $\pi$ is a projective symplectic resolution, $Y$ is a Mori dream space relative to $Z$ \cite[\S3]{AW}. Although the case when $Z$ is affine is treated there, this difference will give no effect on the results except that $\pi$ can have more than one 2-dimensional fibers.

By \cite[Thm. 4.1]{AW}, the classes of $E_1,\dots,E_n$ in the relative Picard group $N^1(Y/Z)$ are linearly independent. Since $\pi$ is the unique symplectic resolution, the (relative) nef cone and the movable cone of $\pi$ are the same. Therefore, the movable cone is strictly convex and the exceptional divisors generate $N^1(Y/Z)$ \cite[Thm. 3.6]{AW}. Thus the claim follows.
\qed

\begin{lem}\label{lem:bundle}
Each $E_i$ is a $\mathbb{P}^1$-bundle over the $K3$ surface $S$ in Lemma \ref{lem:moduli}.
\end{lem}

{\em Proof.} Let $e_i\in N_1(Y/Z)$ be the numerical class of a general fiber of $\pi|_{E_i}$. Then the facets of the (relative) movable cone of $\pi$ are formed by the orthogonal hyperplanes of $e_i$'s with respect to the intersection pairing \cite[Thm. 3.6]{AW}. Since the nef cone coincides with the movable cone, we can choose a $\pi$-nef divisor $D$ such that $(D.e_i)=0$ and $(D.e_j)\ne0$ for $j\ne i$. Then the linear system $|mD|$ for $m\gg0$ gives a contraction $\pi_i:Y\to Z_i$ which contracts the single divisor $E_i$.

As shown in Lemma \ref{lem:moduli}, $Y$ is a Bridgeland moduli space for a $K3$ surface $S$ for some Mukai vector $v$, and therefore any elementary (meaning relative Picard number 1) contractions are realized by ``wall-crossing.'' The complex manifold $\mathrm{Stab}^\dagger(S)$ has wall-and-chamber structure which governs the birational geometry of $Y$. Any wall $\mathcal{W}$ is a codimension one real submanifold of $\mathrm{Stab}^\dagger(S)$ and is associated to a rank-two sublattice $\mathcal{H}$ of $H_\mathrm{alg}^*(S,\alpha,\mathbb{Z}$) (see \cite[\S5]{BM} for details). Elementary divisorial contractions are classified in \cite[Thm. 5.7]{BM} in terms of $\mathcal{H}$. One can check that, in the 4-dimensional case, there are two types depending on whether or not $\mathcal{H}$ contains an isotropic element $w$ (i.e. $w^2=0$) such that $(v.w)=1$ \cite[Lem. 8.7, 8.8]{BM}. In both cases, the exceptional divisor is a $\mathbb{P}^1$-bundle over $S$. 
\qed

\begin{prop}\label{prop:fiber}
Any 2-dimensional fiber $\pi^{-1}(z)$ is isomorphic to a 2-dimensional fiber of a unique symplectic resolution of $\mathrm{Hilb}^2(\Gamma)$ for some $\Gamma$ as in Theorem \ref{thm:main3}.
\end{prop}

{\em Proof.} We first consider the structure of a two-dimensional fiber for $\mathrm{Hilb}^2(\Gamma)$. 

When $\Gamma$ has one singular point, a symplectic resolution of $\mathrm{Hilb}^2(\Gamma)$ are studied and shown to be unique in \cite[\S2]{Y}. Also, the fibers of the 0-dimensional symplectic leaves, which are listed in (\ref{type}), are explicitly described there. (They are isomorphic to some Springer fibers, see \cite[\S4]{Y} and \cite{Lo}.) In any case the fiber consists of irreducible components which are isomorphic to $\mathbb{P}^1\times\mathbb{P}^1$ or the second Hirzebruch surface $\Sigma_2=\mathbb{P}(\mathcal{O}_{\mathbb{P}^1}\oplus \mathcal{O}_{\mathbb{P}^1}(2))$.

When $\Gamma$ has at least two singular points $p_1$ and $p_2$ of type $T_{n_1}$ and $T_{n_2}$ respectively where $T=\mathbf{A},\mathbf{D}$, or $\mathbf{E}$. Then $\mathrm{Hilb}^2(\Gamma)$ has additional 0-dimensional symplectic leaves of the form $\{p_1,p_2\}$ whose analytic type is the same as that of the product singularity $T_{n_1}\times T_{n_2}$. Thus its fiber of the unique symplectic resolution is the product of the two Dynkin trees of projective lines.

We will show that the fiber $F=\pi^{-1}(z)$ is isomorphic to one of the reducible surfaces mentioned above. To do this, we consider the intersection of two $\pi$-exceptional prime divisors $E$ and $E'$. First we assume that the intersection $E\cap E'\cap F$ is 2-dimensional. Since both $E$ and $E'$ are $\mathbb{P}^1$-bundles by Lemma \ref{lem:bundle}, the surface $E\cap E'\cap F$ has two different rulings and thus each connected component $P$ of $E\cap E'\cap F$ is isomorphic to $\mathbb{P}^1\times\mathbb{P}^1$ (As we will see later, $E\cap E'\cap F$ is in fact always connected). One can check that the formal neighborhood of  $P$ in $Y$ is isomorphic to that of the zero section in the cotangent bundle of $\mathbb{P}^1\times\mathbb{P}^1$ \cite[Lem. 3.2]{Y}. Similarly to the argument in the proof of the previous lemma, taking a suitable divisor from the boundary of the movable cone of $\pi$ gives a contraction $\phi$ of $Y$ that contracts exactly $E$ and $E'$. $\phi$ is identified with the natural contraction (affinization)
\begin{equation}\label{P}
T^*(\mathbb{P}^1\times\mathbb{P}^1)\cong T^*\mathbb{P}^1\times T^*\mathbb{P}^1\to(\mathbf{A}_1\text{-sing.})\times(\mathbf{A}_1\text{-sing.})
\end{equation}
in the formal neighborhood of $P$.

Set $\mathrm{Exc}(\pi)=E_1\cup E_2\cup\cdots\cup E_n$. Then a general fiber of $\pi|_{\mathrm{Exc}(\pi)}$ is a Dynkin tree of $\mathbb{P}^1$'s. We will show that a general fiber $C_i$ of $\pi|_{E_i}$ is irreducible, i.e., it is a single $\mathbb{P}^1$ for each $i$. Note that, for general symplectic resolutions, it can happen that a general fiber of a single exceptional divisor is reducible. Such a phenomenon is caused by an automorphism of a Dynkin diagram and of the corresponding general fiber. Possible cases are classified in \cite[Thm. 1.3]{Wi}. Each case is presented as a pair of the type of a Dynkin diagram and a group of automorphisms of the diagram. In our setting, the case $(\mathbf{A}_{2l}, \mathbb{Z}/2\mathbb{Z})$ cannot happen since each $\pi_i$ is a $\mathbb{P}^1$-bundle by the previous lemma. Now we assume that there exists $i$ such that $C_i$ is reducible in order to deduce contradiction. Then $\pi$ is in one of the four cases: $(\mathbf{A}_{2l+1},\mathbb{Z}/2\mathbb{Z})$, $(\mathbf{D}_l,\mathbb{Z}/2\mathbb{Z})$, $(\mathbf{D}_4,\mathfrak{S}_3)$, and $(\mathbf{E}_6,\mathbb{Z}/2\mathbb{Z})$. In any case, there are $i_1$ and $i_2$ such that\\
$\bullet$ $C_{i_1}\cong\mathbb{P}^1$, and $C_{i_2}\cong\mathbb{P}^1\sqcup\mathbb{P}^1$ or $\mathbb{P}^1\sqcup\mathbb{P}^1\sqcup\mathbb{P}^1$\\
$\bullet$ $C_{i_1}$ intersects with each $\mathbb{P}^1$ in $C_{i_2}$ at one point.

Let $\pi_{i_2}:Y\to Z_{i_2}$ be the contraction of $E_{i_2}$ and $\theta:Z_{i_2}\to W$ the contraction of $E_{i_1}$. Such $\theta$ also exists by the same argument in the proof of Lemma \ref{lem:bundle}. Then the restriction of $\theta$ to $\pi_{i_2}(E_{i_2})\cong S$ is generically a covering map $S\to S'$ with the covering transformation group $G\cong\mathbb{Z}/2\mathbb{Z}$ or $\mathbb{Z}/3\mathbb{Z}$ according to the number of the components of $C_{i_2}$.  Note that the $G$-action on an open subset of $S$ extends to the whole of $S$ since a $K3$ surface is a minimal surface. Thus the map $S\to S'$ can be regarded as a quotient map by $G$. Since $E_{i_1}$ is also a $\mathbb{P}^1$-bundle over $S$, the image $S'$ is birational to $S$. We have the following diagram

$$\begin{CD} @. S\\ @. @VV r V\\ S@>q>> S'=S/G\end{CD}$$
where $q$ is the quotient map and $r$ is the minimal resolution. Let $M$ be the sublattice of $H^2(S,\mathbb{Z})$ generated by the exceptional curves of $r$. Then $q$ naturally induces the push-forward map $q_*:H^2(S,\mathbb{Z})\to M^\perp\subset H^2(S,\mathbb{Z})$ and the pullback $q^*:M^\perp\to H^2(S,\mathbb{Z})$ which preserve algebraic cocycles. One sees that the restriction of $q_*$ to the invariant part $H^2(S,\mathbb{Z})^G$ is an isomorphism onto its image while it multiplies the intersection form by $\sharp G$ (cf. \cite[Lem. 3.1]{Mo} or \cite[Thm. 2.1]{Wh}).

Since the quotient of $S$ by $G$ is birational to a $K3$ surface, $G$ preserves the symplectic form by \cite[Ch. 15, Lem. 4.8]{H2}. Therefore, one sees that the transcendental lattice $T_S=NS(S)^\perp \subset H^2(S,\mathbb{Z})$ is inside the invariant part $H^2(S,\mathbb{Z})^G$. In particular, $q_*$ gives an isomorphism $T_S\cong T_S\subset M^\perp$, though it is a nontrivial multiplication map as lattices. This is absurd since the same lattices have different discriminants. Thus we have shown that every $C_i$ is irreducible.

Next, let us take two divisors $E_{i_1},\,E_{i_2}$ and general fibers $C_{i_1},\,C_{i_2}$ of $\pi|_{E_{i_1}}$ and $\pi|_{E_{i_2}}$ respectively such that $C_{i_1}$ and $C_{i_2}$ intersect. Then the restriction of $\theta$ to $\pi_{i_2}(E_{i_2})\cong S$ is birational and hence an isomorphism. This implies that $E_{i_1}\cap E_{i_2}$ is a section of both of the $\mathbb{P}^1$-bundles $\pi|_{E_{i_1}}$ and $\pi|_{E_{i_2}}$. We say that such $E_{i_1}$ and $E_{i_2}$ are {\em adjacent}. If $E_i$ and $E_j$ satisfy $\pi(E_i)=\pi(E_j)$ (or equivalently $E_i$ can reach $E_j$ via adjacent divisors), we say that $E_i$ and $E_j$ are in the same {\em bunch}.

Now we assume that a 2-dimensional fiber $F=\pi^{-1}(z)$ exists. Note that $F$ is connected since $\pi_*\mathcal{O}_Y=\mathcal{O}_Z$. We will determine the structure of $F$. We first assume that the $\pi$-exceptional prime divisors that contain $F$ are in the same bunch. If any two non-adjacent divisors in this bunch do not have 2-dimensional intersection contained in $F$, then the intersection would be empty. Indeed, otherwise the intersection would be a finite union of fibers of the $\mathbb{P}^1$-bundles, which is obviously absurd. In this case, any fiber of $\pi$ in the bunch would be a Dynkin tree of $\mathbb{P}^1$'s, which is a contradiction. Thus we may assume that non-adjacent divisors $E_1$ and $E_2$ contain a 2-dimensional irreducible component $P\subset F$. Then $P$ is isomorphic to $\mathbb{P}^1\times\mathbb{P}^1$ as we have already discussed. Let $E$ be another $\pi$-exceptional prime divisor that is not adjacent to both $E_1$ and $E_2$. If $E$ intersects with $P$, then $E$ must contain $P$ since the intersection numbers of any fibers of $\pi_1$ and $\pi_2$ with $E$ are zero. However, this is a contradiction since $E_i\,(i=1,2)$ and $E$ do not share fibers of their $\mathbb{P}^1$-bundles. On the other hand, if $E$ is adjacent to $E_1$ or $E_2$, then it always intersects with $P$.\\
\\
{\bf Case 1.} $E$ is adjacent to both $E_1$ and $E_2$.\\
In this case the intersection $C=E\cap P$ is a section of $P$ with respect to both $\pi_1|_P$ and $\pi_2|_P$. Therefore, it must be a diagonal (i.e. an irreducible curve with $C^2=2$) in $P$. By using normal bundle sequences and taking the fact that $N_{P/Y}\cong \Omega_P$ into consideration, one sees that $N_{C/Y}\cong \mathcal{O}(2)\oplus \mathcal{O}(-2)^{\oplus 2}$.

On the other hand $E$ contains a $\mathbb{P}^1$-bundle $\Sigma\subset F$ over $C$, and $\Sigma$ contains $C$ as a negative section. By the calculation of the normal bundle above, we see that $\Sigma$ is isomorphic to $\Sigma_2$.\\
\\
{\bf Case 2.} $E$ is adjacent to just one, say $E_1$, of $E_1$ and $E_2$.\\
In this case the intersection $C=E\cap P$ is the section of $\pi_1|_P$. Since the intersection number of any fiber of $\pi_2$ and $E$ is zero, $C$ must be a fiber of $\pi_2|_P$. Thus in this case $N_{C/Y}\cong \mathcal{O}(2)\oplus \mathcal{O}^{\oplus 2}$, and the corresponding $\mathbb{P}^1$-bundle $\Sigma'\subset E\cap F$ over $C$ is isomorphic to $\mathbb{P}^1\times\mathbb{P}^1$. Note that $\Sigma'$ is also contained in $E_2$. Indeed, the rulings of $\Sigma'\cong\mathbb{P}^1\times\mathbb{P}^1$ of the other direction than $E\to S$ are numerically equivalent to $C\subset E_2$. Thus $\Sigma'\subset E_2\cap E\cap F$ is in a similar situation to $P\subset E_1\cap E_2\cap F$.

Next we consider the intersection of $\Sigma$ in {\bf Case 1} and another exceptional divisor. The same argument as above shows that, if a $\pi$-exceptional prime divisor $E'$ other than $E_1,E_2$ and $E$ intersects with $\Sigma$, it is adjacent to one of $E,E_1$ and $E_2$. When $E'$ is adjacent to $E$, such $E'$ is unique because of the shapes of the Dynkin diagrams. Since the intersection number of fibers of $\pi_1$ and $\pi_2$ with $E'$ are zero (otherwise a closed path would appear in a Dynkin diagram), the two sections $C$ and $C':=\Sigma\cap E'$ are disjoint and thus $C'$ is a section of $\Sigma\to\mathbb{P}^1$ with $C'^2=2$. By the normal bundle argument as above shows that the $\mathbb{P}^1$-bundle $\Sigma''\subset E'\cap F$ over $C'$ is isomorphic to $\Sigma_2$. 

When $E'$ is adjacent to $E_1$ or $E_2$, we see that the situation is similar to {\bf Case 2} and $\Sigma$ intersects along its fiber with (a connected component of) $E\cap E' \cap F$ which is isomorphic to $\mathbb{P}^1\times\mathbb{P}^1$.

Let us consider again the case when $E'$ is adjacent to $E$. We can describe the intersections of $\Sigma''$ above with other exceptional divisors in a similar way: any $\pi$-exceptional prime divisor $E''$ other than $E$ and $E'$ which intersects with $\Sigma''$ is adjacent to one of $E',E_1$ and $E_2$. When $E''$ is adjacent to $E'$, the intersection $C''$ is a section which is disjoint from $C'$ with $C''^2=2$. $E''$ contains the $\mathbb{P}^1$-bundle over $C''$ which is isomorphic to $\Sigma_2$ and contains $C''$ as the $(-2)$-curve. When $E''$ is adjacent to $E_i$, the component $\Sigma''$ intersects along its fiber with (a connected component of) $E'\cap E'' \cap F$ which is isomorphic to $\mathbb{P}^1\times\mathbb{P}^1$.
 
Note that every irreducible component of $F$ intersects with at least two exceptional prime divisors since otherwise $F$ could not be a fiber. Since we have considered all intersections with the exceptional divisors starting from $P$, every irreducible component of $F$ will appear in the above procedures. One can check that the resulting $F$, up to isomorphisms, fall into one of the types in the list (\ref{type}) (see \cite[\S2]{Y} or \cite{Lo}). In particular, we see that every component of $F$ which is isomorphic to $\mathbb{P}^1\times\mathbb{P}^1$ is of the form $E_i\cap E_j\cap F$ for some $i,j$.

\begin{figure}[H]
\centering
\includegraphics[scale=0.3]{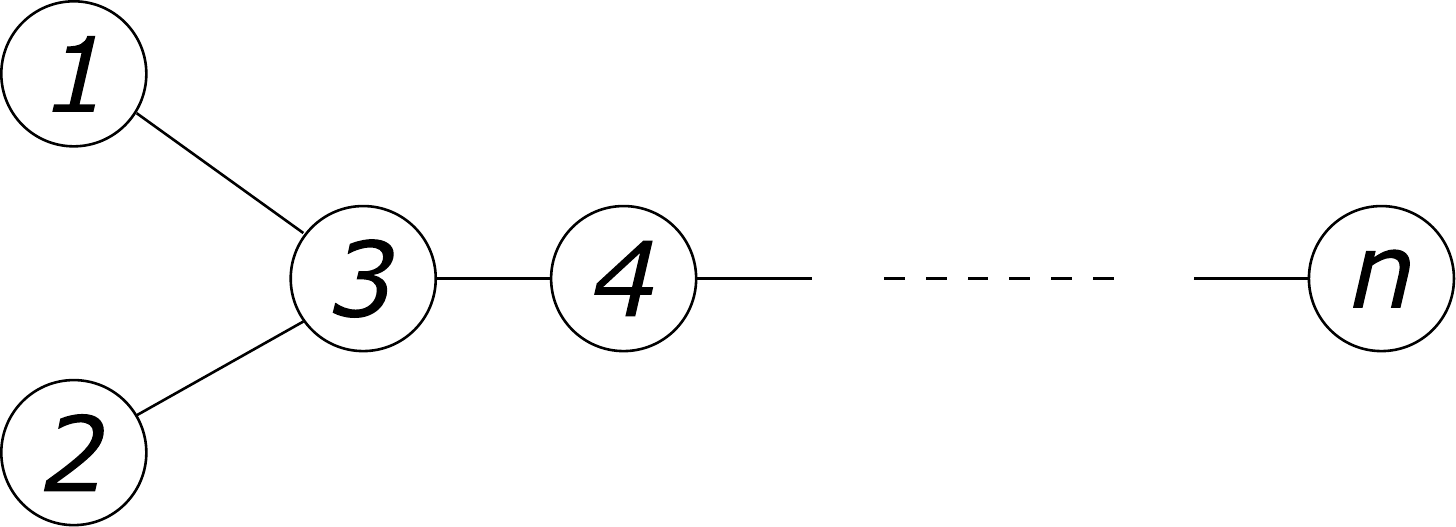}

\caption{The Dynkin diagram of type $\mathbf{D}_n$}
\label{fig1}
\end{figure}

As an example, let us consider the case when a general fiber is of type $\mathbf{D}_n\,(n\ge4)$. Let $E_1,\dots,E_n$ be the irreducible $\pi$-exceptional divisors containing $F$ such that $E_i$ corresponds to the $i$-th vertex of the Dynkin diagram in Figure \ref{fig1}.

We first assume that $E_1\cap E_2$ contains a 2-dimensional component $P$ of $F$. Then $P$ is isomorphic to $\mathbb{P}^1\times\mathbb{P}^1$, and $E_3$ is a unique divisor that intersects with $P$. Also, $E_3$ contains a $\mathbb{P}^1$-bundle $\Sigma$ over a diagonal of $P$ (see {\bf Case 1}). $E_4$ is a unique divisor that intersects with $\Sigma$ and contains $\mathbb{P}^1$-bundle $\Sigma'$ over a curve $C$ of $\Sigma$ with $C^2=2$. After all, we will obtain a chain of one $\mathbb{P}^1\times\mathbb{P}^1$ and $n-2$ copies of $\Sigma_2$. This corresponds to $\bar{\mathbf{D}}_{n,\mathrm{I}}$.

Next we consider the case when $E_1\cap E_2$ does not contain a 2-dimensional component of $F$. When $n=4$, we can reduce to the previous case by symmetry. We assume $n\ge5$. In this case one can see that every pair of non-adjacent divisors intersect along a 2-dimensional component of $F$. Indeed, for any pair of such 2-dimensional components $P_{i,j}\subset E_i\cap E_j$, there is a chain of $P_{i,j}$'s that connects the pair. For example, $P_{1,4}$ and $P_{2,4}$ are connected by the chain 
$P_{1,4},\,P_{1,5},\,P_{3,5},\,P_{2,5},\,P_{2,4}$. Finally we will find that $F$ corresponds to $\bar{\mathbf{D}}_{n,\mathrm{I\hspace{-1pt}I}}$.

To finish the classification, let us consider the case when $F$ is contained in more than one bunches of exceptional divisors. In this case there are two divisors, say $E_1$ and $E_2$, that are not in the same bunch and $E_1\cap E_2$ contains a 2-dimensional component $P$ of $F$. $P$ must be isomorphic to $\mathbb{P}^1\times\mathbb{P}^1$ as discussed. Note that the number of bunches cannot be more than two. This is shown as follows. Assume there is an exceptional divisor $E_3$ in a third bunch. Since the intersection numbers of $E_3$ and general fibers of $E_1$ and $E_2$ are zero, $P$ is contained in $E_3$. However, this is contrary to the fact that, in the neighborhood of $P$, exactly two divisors are contracted (cf. (\ref{P})).

We can use the same method as above to produce the other components of $F$. In this case all components of $F$ are isomorphic to $\mathbb{P}^1\times\mathbb{P}^1$ since no divisors are adjacent to both $E_1$ and $E_2$ and thus {\bf Case 1} will not happen. It is easy to see that the resulting $F$ is isomorphic to the product of two Dynkin trees of $\mathbb{P}^1$. This completes the proof.
\qed

\vspace{3mm}

\noindent{\em Proof of Theorem \ref{thm:main3}}

By Proposition \ref{prop:fiber}, we have shown that $(Z,z)$ and $(\mathrm{Hilb}^2(\Gamma),q)$ have isomorphic 2-dimensional fibers in their unique symplectic resolutions for some $q$. When the fiber is of one of the types in the list (\ref{type}), this implies that $(Z,z)$ and $(\mathrm{Hilb}^2(\Gamma),q)$ have the same singularity types \cite[Thm. 3.4]{Y}.

In the product cases, we will prove the claim using the results in \cite[\S3]{Y}. Every irreducible component of $F\subset Y$ is isomorphic to $\mathbb{P}^1\times\mathbb{P}^1$ and thus its formal neighborhood is always isomorphic to the completion of the cotangent bundle along the zero section. We will show that the isomorphism class of the formal neighborhood of $F$ is also determined uniquely.

As we have seen, the fiber $F\subset Y$ is contained in two bunches of exceptional divisors. Let $B_1=E_1\cup\cdots\cup E_k$ and $B_2=E'_1\cup\dots\cup E'_l$ be the two bunches. There are contractions $\pi_{B_1}$ and $\pi_{B_2}$ of $B_1$ and $B_2$ respectively. Then the images $S_1:=\pi_{B_1}(B_1)$ and $S_2:=\pi_{B_2}$ are isomorphic to the minimal resolutions of ADE-singularities. Let $T_l$ (resp. $T_k$) be the ADE-type of $S_1$ (resp. $S_2$). The unique symplectic resolution of the product of two ADE-singularities of types $T_l$ and $T_k$ also admits two bunches of divisors and corresponding surfaces $S'_1$ and $S'_2$ which are isomorphic as symplectic varieties to $S_1$ and $S_2$ respectively in the formal neighborhoods of the exceptional curves.

We fix the two symplectic isomorphisms between $S_1\cong S'_1$ and $S_2\cong S'_2$. Then, by applying \cite[Lem. 3.5 and Rem. 3.4]{Y} repeatedly, we obtain a (non-canonical) isomorphism of the formal neighborhoods of central fibers of $Y$ and the product resolution which induces the given isomorphisms of $S_i$'s and $S'_i$'s. This implies that the germ $(Z,z)$ is isomorphic to the product of two ADE-singularities. As we already see in the proof of Proposition \ref{prop:fiber}, such a germ is also presented as $(\mathrm{Hilb}^2(\Gamma),q)$ for some $\Gamma$ and $q$. 
\qed

\vspace{3mm}

{\bf Remark.}  As we will see below, Theorem \ref{thm:main3} can apply to $F(X)$. The proof of this theorem shows that three or more irreducible components of the singular locus $\mathrm{Sing}(Z)$ of $Z$ do not intersect at one point. This was discussed for $F(X)$ in Lemma \ref{S_p}. Moreover, from the proof we also know that all irreducible components of $\mathrm{Sing}(Z)$ are birational to each other. This is also a generalization of the result for $F(X)$ \cite[Prop. 5.8(c)]{O'G}.

\vspace{3mm}

\noindent{\em Proof of Theorem \ref{thm:main1} and \ref{thm:main2}}

As mentioned in the previous section, we have a birational morphism $\mathrm{Hilb}^2(S_p)\to F(X)$. Since $S_p$ is a surface with only ADE-singularities, $\mathrm{Hilb}^2(S_p)$ has a symplectic resolution $Y$. By \cite[Thm. 27.8]{GHJ}, we see that $Y$ is of $K3^{[2]}$-type.

To show that $Y\to F(X)$ is a unique symplectic resolution, it suffices to show that two-dimensional fibers do not contain $\mathbb{P}^2$ since two different symplectic resolutions are related by a sequence of Mukai flops \cite[Thm. 1.2]{WW}. We already know that the fibers of $Y\to\mathrm{Hilb}^2(S_p)$ do not contain $\mathbb{P}^2$. One can check that the fiber of $\mathrm{Hilb}^2(S_p)\to F(X)$ of a point $x\in S_p\subset F(X)$ can be identified with the set of lines passing through $x$ that are contained in the quadric $\mathcal{Q}$. Thus it is 1-dimensional or isomorphic to a non-degenerate or corank-1 quadric in $\mathbb{P}^3$. Therefore, no irreducible components which is isomorphic to $\mathbb{P}^2$ cannot appear in the fiber.

Now we can apply Theorem \ref{thm:main3} to $Y\to F(X)$. Let
$x\in F(X)$ be a 0-dimensional symplectic leaf. If $\mathrm{Sing}(F(X))$ is reducible at $x$, then we fall into the product case. The types of $\Gamma_i$'s are determined by the types of $x$ in $S_{p_i}$'s, which are $T_{n_i}$'s. This proves Theorem \ref{thm:main1}.

If $\mathrm{Sing}(F(X))$ is irreducible at $x$, then we fall into one of the types in the list (\ref{type}). In this case, however, the types in (\ref{type}) are not made distinct just by looking at the type of $x$ in $S_p$. For example, $\bar{\mathbf{A}}_7$ and $\bar{\mathbf{E}}_6$ both have $\mathbf{A}_5$ as the types in the 2-dimensional symplectic leaves. By taking the number of the irreducible exceptional divisors of $Y\to F(X)$ also into account, we can uniquely determine the type of $x\in F(X)$. Note that the resolution $Y\to F(X)$ factors through $\mathrm{Hilb}^2(S_p)$, and thus we see that the number of the irreducible exceptional divisors is equal to $n$ for type $T_n$. This proves Theorem \ref{thm:main2}.
\qed

\vspace{0.2cm}

\begin{center}
Department of Mathematics, Graduate School of Science, Kyoto University, Japan 

ryo-yama@math.kyoto-u.ac.jp
\end{center}

\end{document}